\documentclass[10pt,oneside,a4paper]{amsart} 
\usepackage{latexsym}
\usepackage{amsfonts,amsmath,amssymb,indentfirst}
\usepackage[english]{babel}
\newcommand{\beq}{\begin{equation}}
\newcommand{\eeq}{\end{equation}}
\newcommand{\bdism}{\begin{displaymath}}
\newcommand{\edism}{\end{displaymath}}

\newcommand{\setS}{\mathbb S}
\newtheorem{theorem}{Theorem}[section]
\newtheorem{proposition}[theorem]{Proposition}
\newtheorem{corollary}[theorem]{Corollary}

\newtheorem{remark}[theorem]{Remark}
\newtheorem{conjecture}[theorem]{Conjecture}

\author{\scshape Luca Fabrizio Di Cerbo*}
\title{\bf Eigenvalues of the Laplacian under the Ricci Flow}
\begin{document}
\pagestyle{headings} 
\thanks{*Supported in part by a Renaissance Technologies Fellowship.}
\address{Department of Mathematics, SUNY, Stony Brook, NY 11794-3651,
USA} \email{luca@math.sunysb.edu}
\begin{abstract}
We derive, under a technical assumption, the first variation
formula for the eigenvalues of the Laplacian on a closed manifold
evolving by the Ricci flow and give some applications.
\end{abstract}
\maketitle

\section{Introduction}
\pagenumbering{arabic}

Given a closed manifold $M^{n}$ endowed with a Riemannian metric
$g_{0}$, the Ricci flow determines a one parameter family of
metrics $g(t)$ via the geometric evolution equation
\begin{align}\label{Ricci equation}
\frac{\partial g}{\partial t}=-2Ric(g)
\end{align}
with initial condition \bdism g(0)=g_{0}. \edism Short time
existence and uniqueness for solution to the Ricci flow was first
shown in Hamilton's seminal paper \cite{Ha1}, using the Nash-Moser
theorem, and shortly after by DeTurck in \cite{De} who
substantially simplified the proof. Strictly related to the Ricci
flow is the evolution equation
\begin{align}\label{Ricci norm}
\frac{\partial g}{\partial t}=\frac{2}{n}rg-2Ric(g)
\end{align}
where $r=\frac{\int_{M}Rd\mu}{\int_{M}d\mu}$, which is called
normalized Ricci flow since has the remarkable property to
preserve the volume of the initial Riemannian manifold. As shown
in \cite{Ha1}, the evolution equations \ref{Ricci equation} and
\ref{Ricci norm} differ only by a change of scale in space and a
change of parametrization in time. Since its introduction, the
Ricci flow has been a very effective tool for studying the
topology of manifolds. In the most of first applications
\cite{Ha1}, \cite{Ha2}, \cite{Ha3}, \cite{Hu}, the normalized flow
is shown to converge to a canonical metric, e.g. in \cite{Ha1}
Hamilton proves that the normalized flow on any closed
$3$-manifold of positive Ricci curvature exists for all time and
converges to a positive constant sectional curvature metric. In
this sense, the Ricci flow may be regarded as a natural homotopy
between a given metric of positive Ricci curvature and a canonical
metric of constant sectional curvature in the same volume class.
Analogously, whenever the Ricci flow converges, it can be
considered as a natural homotopy between the initial metric and
the limit metric that is necessarily Einstein \bdism
Ric(g_{\infty})=\frac{1}{n}\frac{\int_{M}R_{\infty}d\mu_{\infty}}{\int_{M}d\mu_{\infty}}g_{\infty}.
\edism It is then natural to use the Ricci flow to study
properties of geometrical meaningful objects, such as the
eigenvalues of the Laplacian. In this paper we study, under a
technical assumption, the behavior of the spectrum of the
Laplacian when we deform a given initial metric by the
Ricci flow.\\
Throughout the paper, let $M^{n}$ be a closed manifold of
dimension $n$, and $C^{\infty}(S^{+}_{2}TM)$ the space of smooth
Riemannian metrics on $M^{n}$. For $g\in C^{\infty}(S^{+}_{2}TM)$,
let $\Delta_{g}=Tr_{g}\nabla^{2}$ be the Laplacian operator on
$C^{\infty}(M)$ and \bdism
Spec(g)=\left\{0=\lambda_{0}(g)<\lambda_{1}(g)\leq\lambda_{2}(g)\leq...\leq\lambda_{k}(g)\leq...\right\}
\edism
the spectrum of $\Delta_{g}$.\\
Given $g_{0}\in C^{\infty}(S^{+}_{2}TM)$, let $g(t)$ be the smooth
$1$-parameter family of metrics who solves the Ricci flow with
initial metric $g_{0}$. We can then regard each eigenvalue
$\lambda(g(t))$ as a function of one real parameter, and try to
study its evolution at least locally in time. This problem is
clearly related to the general problem concerning the behavior of
the spectrum under a suitable deformation $g(t)$ of $g_{0}$ in
$C^{\infty}(S^{+}_{2}TM)$. In this direction a classical result is
the following
\begin{theorem}
For $g\in C^{\infty}(S^{+}_{2}TM)$ and $h\in C^{\infty}(S_{2}TM)$,
let $g(t)=g+th$, $\left|t\right|<\epsilon$ for sufficiently small
$\epsilon>0$. Let $\lambda$ be an eigenvalue of $\Delta_{g}$ with
multiplicity $l$. Then there exist $f_{i}(t)\in C^{\infty}(M)$,
i=1,...,l, such that
\begin{enumerate}
\item $\lambda_{i}(t)$ and $f_{i}(t)$ depend real analytically on t,
$\left|t\right|<\epsilon$, for each $i=1,...,l$,
\item $\Delta_{g(t)}f_{i}(t)+\lambda_{i}(t)f_{i}(t)=0$, for each
$i=1,...,l$ and $t$,
\item $\lambda_{i}(0)=\lambda$, $i=1,...,l$, and
\item $\left\{f_{i}(t)\right\}^{l}_{i=1}$ is orthonormal with respect
to $\left\langle, \right\rangle_{g(t)}$ for each t.
\end{enumerate}
\end{theorem}
The above theorem is due to Berger \cite{Berg1}, and it is one of
the first attempt in studying the behavior of the spectrum under a
perturbation of the metric. On the other hand Bando and Urakawa
have extended the above theorem to a family of metrics which
depends real analytically on time, see \cite{Ura1}. Unfortunately,
the theorems of Berger and Bando-Urakawa does not apply directly
to a smooth $1$-parameter family of metrics generated by solving
an IVP of the Ricci flow type. Since the Ricci flow equation is a
parabolic equation for the metric see \cite{Ha1}, its solutions
are expected to not depend real analytically on the time variable.
More precisely, the Ricci flow equation implies that the scalar
curvature function of the evolving metric satisfies an heat-type
scalar equation. Hence if the solution flow were analytic in time,
we could analytically continue backwards in time, which implies
that the associated scalar curvature
equation admits solution for negative time.\\
In what follows we then assume that, under a Ricci flow
deformation $g(t)$ of a given initial metric, a result of the
Berger-Bando-Urakawa type holds, that is we assume the existence
and $C^{1}$-differentibility
of the elements $\lambda_{i}(t)$ and $f_{i}(t)$.\\
We remark that the problem of the spectrum variation under a
deformation of the metric given by a parabolic equation seems not
to be presented explicitely before. However in the celebrated
paper \cite{Perelman}, a similar geometric problem is cosidered
from a different point of view. In this paper Perelman introduces
a riemannian functional whose gradient flow is the Ricci flow
modulo diffeomorphisms. It turns out that this functional is the
lowest eigenvalue $\lambda_{1}$ of the operator $-4\Delta+R$,
therefore proving that the quantity $\lambda_{1}(g(t))$ is
nondecreasing along the flow. For the details and applications of
this important observation we refer to the original paper.

\section{Variation Formulas}
Let $(M^{n},g(t))$ be a solution of the Ricci flow on the maximal
time interval $[0,T)$ and consider the associated Rayleigh-Ritz
quotient \bdism \lambda(t)=\frac{\int_{M}{\left|\nabla
f_{t}\right|}^{2}d\mu_{t}}{\int_{M}f^{2}_{t}d\mu_{t}}, \edism
which defines the evolution of an eigenvalue of the Laplacian
under the flow. We can then consider the change rate of any
eigenvalue taking the time derivative of the Rayleigh-Ritz
quotient which defines it. For simplicity we consider normalized
eigenfunctions i.e. \beq\label{marchia}
\int_{M}f_{t}d\mu=0,\\
\int_{M}f^{2}_{t}d\mu=1, \eeq then
\begin{align} \notag
\frac{d\lambda}{dt}=&\int_{M}\frac{d}{dt}{\left|\nabla
f\right|}^{2}d\mu+\int_{M}{\left|\nabla
f\right|}^{2}\frac{d}{dt}d\mu \\ \notag
=&\int_{M}\frac{d}{dt}(g^{ij})\nabla_{i}f\nabla_{j}fd\mu+2\int_{M}\left\langle
\nabla f^{'},\nabla f\right\rangle d\mu-\int_{M}R{\left|\nabla
f\right|}^{2}d\mu \\ \notag
=&2\int_{M}g^{ik}g^{jl}R_{kl}\nabla_{i}f\nabla_{j}fd\mu+2\int_{M}\left\langle
\nabla f^{'},\nabla f\right\rangle d\mu-\int_{M}R{\left|\nabla
f\right|}^{2}d\mu \\ \notag =&2\int_{M}Ric(\nabla f,\nabla
f)d\mu+2\int_{M}\left\langle \nabla f^{'},\nabla f\right\rangle
d\mu-\int_{M}R{\left|\nabla f\right|}^{2}d\mu,
\end{align}
for the time derivative of the volume element see \cite{Ha1}. Now,
using \ref{marchia} we get the following two integrability
conditions
\begin{align}
&\int_{M}f^{'}d\mu=\int_{M}fRd\mu \label{integrability1}\\
&\int_{M}f^{'}fd\mu=\frac{1}{2}\int_{M}f^{2}Rd\mu\label{integrability2}.
\end{align}
Finally, by \ref{integrability2} we have
\begin{align} \notag
2\int_{M}\left\langle \nabla f^{'},\nabla f\right\rangle
d\mu=&-2\int_{M}f^{'}\Delta f d\mu \\ \notag =&2\lambda
\int_{M}f^{'}f d\mu=\lambda \int_{M}f^{2}R d\mu. \\ \notag
\end{align}
We have thus proved the following proposition
\begin{proposition}\label{evolution1}
Let $(M^{n},g(t))$ be a solution of the unnormalized Ricci on the
smooth manifold $(M^{n},g_{0})$. If $\lambda(t)$ denotes the
evolution of an eigenvalue under the Ricci flow, then \bdism
\frac{d\lambda}{dt}=\lambda \int_{M}f^{2}R
d\mu-\int_{M}R{\left|\nabla f\right|}^{2}d\mu+2\int_{M}Ric(\nabla
f,\nabla f)d\mu. \edism where $f$ is the associated normalized
evolving eigenfunction.
\end{proposition}
As shown below, an analogous equation holds for the normalized
Ricci flow
\begin{proposition}\label{evolution2}
Let $(M^{n},g(t))$ be a solution of the normalized Ricci on the
smooth manifold $(M^{n},g_{0})$. If $\lambda(t)$ denotes the
evolution of an eigenvalue under the normalized Ricci flow, then
\bdism \frac{d\lambda}{dt}=-\frac{2}{n}r\lambda+\lambda
\int_{M}f^{2}R d\mu-\int_{M}R{\left|\nabla
f\right|}^{2}d\mu+2\int_{M}Ric(\nabla f,\nabla f)d\mu. \edism
where $f$ is the associated normalized evolving eigenfunction.\\
\end{proposition}
\begin{proof}
In the normalized case, the integrability conditions read as
follows
\begin{align} \label{bella1}
&\int_{M}f^{'}d\mu=\int_{M}fRd\mu, \\ \label{bella2}
&2\int_{M}ff^{'}d\mu=\int_{M}f^{2}Rd\mu-r,
\end{align}
since \bdism
\frac{d}{dt}d\mu=\frac{1}{2}Tr_{g}(\frac{2}{n}rg-2Ric)d\mu=(r-R)d\mu.
\edism We can then write
\begin{align}\notag
\frac{d\lambda}{dt}=&\int_{M}\frac{d}{dt}{\left|\nabla
f\right|}^{2}d\mu+\int_{M}{\left|\nabla
f\right|}^{2}\frac{d}{dt}d\mu \\ \notag
=&\int_{M}\frac{d}{dt}(g^{ij})\nabla_{i}f\nabla_{j}fd\mu+2\int_{M}\left\langle
\nabla f^{'},\nabla f\right\rangle d\mu\\ \notag
&+\int_{M}(r-R){\left|\nabla f\right|}^{2}d\mu \\ \notag
=&-\int_{M}g^{ik}g^{jl}(\frac{2r}{n}g_{kl}-2R_{kl})\nabla_{i}f\nabla_{j}fd\mu+2\int_{M}\left\langle
\nabla f^{'},\nabla f\right\rangle d\mu\\ \notag
&+\int_{M}(r-R){\left|\nabla f\right|}^{2}d\mu \\ \notag
=&2\int_{M}Ric(\nabla f,\nabla
f)d\mu-\frac{2r}{n}\int_{M}\left|\nabla
f\right|^{2}d\mu+2\lambda\int_{M}ff^{'}d\mu\\ \notag
&+\int_{M}(r-R){\left|\nabla f\right|}^{2}d\mu \\ \notag
=&2\int_{M}Ric(\nabla f,\nabla
f)d\mu-\frac{2r}{n}\lambda+\lambda\left(\int_{M}f^{2}Rd\mu-r\right)+r\lambda\\
\notag &-\int_{M}R\left|\nabla f\right|^{2}d\mu \\ \notag
=&-\frac{2}{n}r\lambda+\lambda \int_{M}f^{2}R
d\mu-\int_{M}R{\left|\nabla f\right|}^{2}d\mu +2\int_{M}Ric(\nabla
f,\nabla f)d\mu.
\end{align}
\end{proof}

It is now interesting to write down proposition \ref{evolution1}
and proposition \ref{evolution2} in some remarkable particular
case.
\begin{corollary}\label{surface}
Let $(M^{2},g(t))$ be a solution of the unnormalized Ricci flow on
a closed surface, then \bdism
\frac{d\lambda}{dt}=\lambda\int_{M}f^{2}Rd\mu. \edism
\end{corollary}
\begin{proof}
In dimension $n=2$ we have \bdism Ric=\frac{1}{2}Rg, \edism then
\begin{align}\notag
\frac{d\lambda}{dt}&=\lambda \int_{M}f^{2}R
d\mu-\int_{M}R{\left|\nabla f\right|}^{2}d\mu+2\int_{M}Ric(\nabla
f,\nabla f)d\mu\\ \notag
&=\lambda\int_{M}f^{2}Rd\mu-\int_{M}R{\left|\nabla
f\right|}^{2}d\mu+\int_{M}R{\left|\nabla f\right|}^{2}d\mu \\
\notag &=\lambda\int_{M}f^{2}Rd\mu.
\end{align}
\end{proof}
\begin{corollary}\label{variation1}
Let $(M^{2},g(t))$ be a solution of the normalized Ricci flow on a
closed surface with normalized initial metric, then
\begin{align}\notag
\frac{d\lambda}{dt}&=\lambda\int_{M}f^{2}Rd\mu-r\lambda=\lambda\int_{M}R(f^{2}-1)d\mu
\\ \notag
\end{align}
\end{corollary}
\begin{proof}
Obvious.
\end{proof}
\begin{remark}
Because of the Gauss-Bonnet theorem the above variation formula
can be written as \bdism
\frac{d\lambda}{dt}=\lambda\int_{M}f^{2}Rd\mu-\lambda4\pi\chi(M),
\edism where $\chi(M)$ is the Euler characteristic of the surface.
\end{remark}

Let us now consider the behavior of the spectrum when we evolve an
initial metric that is homogeneous. As shown in \cite{Ha1}, the
Ricci flow preserves the isometries of the initial Riemannian
manifold. We conclude that the evolving metric remains homogeneous
during the flow. This important observation implies the following
\begin{corollary}\label{homogeneous1}
Let $(M^{n},g(t))$ be a solution of the unnormalized Ricci on the
smooth homogeneous manifold $(M^{n},g_{0})$. If $\lambda(t)$
denotes the evolution of an eigenvalue under the Ricci flow, then
\bdism \frac{d\lambda}{dt}=2\int_{M}Ric(\nabla f,\nabla f)d\mu.
\edism
\end{corollary}
\begin{proof}
Since the evolving metric remains homogeneous, the thesis follows
from proposition \ref{evolution1} and the fact that an homogeneous
manifold has constant scalar curvature.
\end{proof}
\begin{corollary}
Let $(M^{n},g(t))$ be a solution of the normalized Ricci on the
smooth homogeneous manifold $(M^{n},g_{0})$. If $\lambda(t)$
denotes the evolution of an eigenvalue under the Ricci flow, then
\bdism
\frac{d\lambda}{dt}=-\frac{2}{n}R\lambda+2\int_{M}Ric(\nabla
f,\nabla f)d\mu. \edism
\end{corollary}
\begin{proof}
In the homogeneous case the average scalar curvature can be
written as \bdism
r=\frac{\int_{M}R d\mu}{\int_{M}d\mu}=R.\\
\edism
\end{proof}

\section{Applications}
In this section we show how the variational formulas can be
effectively applied to derive some interesting properties of the
evolving spectrum. First we concentrate on $3$-manifolds, next we
discuss the Riemannian
surface and in particular the $2$-sphere.\\
Let $(M^{3},g_{0})$ be a closed three manifold with positive Ricci
curvature. It is well known that the Ricci flow on such manifold
exists on a limited maximal time interval $[0,T)$, see \cite{Ha1};
we shall then show that the eigenvalues of the Laplacian diverges
as $t\rightarrow T$. The result is independent of the technical
assumption explained in the
introduction.\\
\begin{proposition}
Let $(M^{3},g(t))$ be a solution of the Ricci flow on a closed
$3$-manifold whose Ricci curvature is positive initially, then
\bdism lim_{t\rightarrow T}\lambda(t)=\infty. \edism
\end{proposition}
\begin{proof}
On a closed manifold $M^{n}$, for any smooth functions $f$ holds
the celebrated Reilly formula \bdism \int_{M}\left|\nabla\nabla
f\right|^{2}d\mu+\int_{M}Ric(\nabla f,\nabla
f)d\mu=\int_{M}(\Delta f)^{2}d\mu. \edism Since \bdism
\left|\nabla\nabla f\right|^{2}\geq\frac{1}{n}(\Delta f)^{2},
\edism we have the inequality \beq \label{basta}
\frac{n-1}{n}\int_{M}(\Delta f)^{2}d\mu\geq\int_{M}Ric(\nabla
f,\nabla f)d\mu. \eeq For any solution of the Ricci flow on a
closed three manifold with positive Ricci curvature there exists
$\epsilon>0$ such that the condition \bdism Ric\geq\epsilon Rg
\edism is preserved along the flow, see \cite{Ha1}. Therefore
\bdism \frac{2}{3}\lambda^{2}(t)\geq\int_{M}Ric(\nabla f,\nabla
f)d\mu\geq\epsilon\int_{M}R\left|\nabla
f\right|^{2}d\mu\geq\epsilon R_{min}(t)\lambda(t), \edism and then
\bdism \lambda(t)\geq\frac{3}{2}\epsilon R_{min}(t). \edism The
thesis follows since \bdism lim_{t\rightarrow T}R_{min}(t)=\infty,
\edism see again \cite{Ha1}.
\end{proof}

\begin{proposition}
Let $(M^{3},g(t))$ be a solution to the Ricci flow on a closed
manifold whose Ricci curvature is positive initially. Then there
exists $\overline{t}\in [0,T)$ depending on $g_{0}$ such that for
each $t\in[\overline{t},T)$ the eigenvalues of the Laplacian are
increasing.
\end{proposition}
\begin{proof}
Let $0<\epsilon\leq\frac{1}{3}$ be a constant such that \bdism
Ric\geq\epsilon Rg \edism is preserved under the flow. By
proposition \ref{evolution1} we can then write
\begin{align}\notag
\frac{d\lambda}{dt}&\geq\lambda \int_{M}f^{2}R
d\mu-\int_{M}R{\left|\nabla
f\right|}^{2}d\mu+2\epsilon\int_{M}R\left|\nabla f\right|^{2}d\mu
\\ \notag
&\geq\lambda\left\{R_{min}(t)+(2\epsilon-1)R_{max}(t)\right\}.
\end{align}
As proved in \cite{Ha1}, for each $\eta>0$, we can find
$T_{\eta}\in[0,T)$ such that for $t\in [T_{\eta},T)$ \bdism R\geq
(1-\eta)R_{max} \edism hence the proposition follows easily.
\end{proof}

\begin{proposition}
Let $(M^{3},g(t))$ be a solution of the Ricci flow on a closed
homogeneous $3$-manifold whose Ricci curvature is nonnegative
initially, then the eigenvalues of the Laplacian are increasing.
\end{proposition}
\begin{proof}
Since in dimension three the nonnegativity of the Ricci tensor is
preserved under the Ricci flow, the thesis follows easily from
corollary \ref{homogeneous1}.
\end{proof}
\begin{proposition}
Let $(M^{n},g(t))$ be a solution of the Ricci flow on a closed
homogeneous $n$-manifold whose curvature operator is nonnegative
initially, then the eigenvalues of the Laplacian are increasing.
\end{proposition}
\begin{proof}
As shown in \cite{Ha2}, in any dimension the nonnegativity of the
curvature operator is preserved along the Ricci flow. Since the
nonnegativity of the curvature operator implies the nonnegativity
of the Ricci tensor the thesis follows from corollary
\ref{homogeneous1}.
\end{proof}

\begin{proposition}
Let $(M^{n},g(t))$ be a solution of the Ricci flow on a closed
homogeneous manifold, then \bdism
\frac{d\lambda}{dt}\leq2\frac{(n-1)}{n}\lambda^{2}. \edism
\end{proposition}
\begin{proof}
By \ref{basta} \bdism \int_{M}Ric(\nabla f,\nabla
f)d\mu\leq\frac{n-1}{n}\int_{M}(\Delta f)^{2}, \edism we have thus
\begin{align}\notag
\frac{d\lambda}{dt}&=2\int_{M}Ric(\nabla f,\nabla
f)d\mu\leq2\frac{(n-1)}{n}\int_{M}(\Delta f)^{2} \\ \notag
&=2\frac{(n-1)}{n}\lambda^{2}.
\end{align}
\end{proof}

We conclude studying the spectrum of Riemannian surface evolving
under the flow. Our result is the following
\begin{proposition}\label{monotone}
Let $(M^{2},g_{0})$ be a closed surface with nonnegative scalar
curvature, then the eigenvalues of the Laplacian are increasing
under the Ricci flow.
\end{proposition}
\begin{proof}
Under the unnormalized Ricci flow on a surface, we have \bdism
\frac{\partial}{\partial t}R=\Delta R+R^{2}, \edism see
\cite{Chow1}. By the scalar maximum principle, the nonnegativity
of the scalar curvature is preserved along the flow. We conclude
using corollary \ref{surface}.
\end{proof}
Let us now consider the normalized flow on a surface and the
associated variation of the eigenvalues of the Laplacian. Let
$(\setS^{2},g)$ be a closed Riemannian manifold diffeomorphic to
the two dimensional sphere. In \cite{Chow2}, Chow proved that the
normalized Ricci flow on $(\setS^{2},g)$ exists for all time and
converges to a smooth metric $g_{S}$ of constant curvature, that
is the standard metric. On the other hand, a classical result of
Hersch \cite{Hersch} states that $g_{S}$ maximizes $\lambda_{1}$
among all the Riemannian metrics of the same volume. In
particular, if $g_{S}$ is the standard metric on $\setS^{2}$
normalized to volume one, we have \bdism \lambda_{1}(g_{S})=8\pi,
\edism see \cite{Gal}. We then conjecture that the smallest
positive eigenvalue of the Laplacian $\lambda_{1}(t)$ is
monotonically increasing along the Ricci flow $g(t)$ on
$\setS^{2}$.
\begin{conjecture}
Let $(\setS^{2},g_{0})$ be a topological sphere endowed with a
smooth metric normalized to volume one and let $g(t)$ be the
unique solution of the normalized Ricci flow
\begin{align}\notag
&\frac{\partial g}{\partial t}=(r-R)g \\ \notag &g(0)=g_{0}.
\end{align}
Then $\lambda_{1}(g(t))$ is increasing for all $t\in [0,\infty)$
and converges to $\lambda_{1}(g_{S})=8\pi$.
\end{conjecture}
 Actually the convergence of $\lambda_{1}(g(t))$ to
$\lambda_{1}(g_{S})=8\pi$ as $t\rightarrow\infty$ follows easily
from the fact that $g(t)\rightarrow g_{S}$ in the
$C^{\infty}$-topology, see \cite{Chow2}. In fact, it is simple to
prove that the spectrum of the Laplacian depends continuously on
the metric even with respect to the $C^{0}$-topology, see
again \cite{Gal}.\\
 Unfortunately, there are some difficulties in proving the above
conjecture. It is in fact well known that the standard metric on
$\setS^{2}$ does not maximize $\lambda_{2}$, see \cite{Nad}. Hence
it is not possible to prove in general that the variation formula
of corollary \ref{variation1} is positive on $\setS^{2}$. This
suggests that the "shape" of the evolving first eigenfunction of
the Laplacian plays a fundamental role in proving such kind of
conjecture.

\section{examples}
In this section we determine the behavior of the evolving spectrum
on self-similar solutions to the Ricci flow, which are called
\emph{Ricci soliton}. Let $(M^{n},g(t))$ be a solution to the
Ricci flow with initial condition $g(0)=g_{0}$. The solution
$g(t)$ is called Ricci soliton if there exist a smooth function
$\sigma(t)$ and a $1$-parameter family of diffeomorphisms
$\left\{\psi_{t}\right\}$ of $M^{n}$ such that
\begin{align}\label{ricci soliton}
g(t)=\sigma(t)\psi^{*}_{t}(g_{0})
\end{align}
with $\sigma(0)=1$ and $\psi_{0}=id_{M^{n}}$. Now, taking the time
derivative in \ref{ricci soliton} and evaluating the result for
$t=0$, we get that the metric $g_{0}$ satisfies the identity
\begin{align}\label{lie soliton}
-2Ric(g_{0})=2\epsilon g_{0}+L_{X}g_{0},
\end{align}
where $\epsilon=\frac{\sigma^{'}(0)}{2}$ and $X$ is the vector
field on $M^{n}$ generated by $\left\{\psi_{t}\right\}$ for $t=0$.
Conversely, given a metric $g_{0}$ which satisfies \ref{lie
soliton} there exist $\sigma(t)$ and $\left\{\psi_{t}\right\}$
such that $g(t)=\sigma(t)\psi^{*}_{t}(g_{0})$ is a solution to the
Ricci flow with initial condition $g(0)=g_{0}$. In particular we
can choose $\sigma(t)=1+2\epsilon t$ and $\left\{\psi_{t}\right\}$
as the $1$-parameter family of diffeomorphisms generated by the
vector fields \bdism Y_{t}=\frac{1}{\sigma(t)}X, \edism see
\cite{Chow1}. In summary, each self-similar solution to the Ricci
flow can be written in the canonical form \bdism g(t)=(1+2\epsilon
t)\psi^{*}_{t}(g_{0}). \edism We then say that the soliton is
expanding, shrinking, or steady, if
$\epsilon>0$, $\epsilon<0$, or $\epsilon=0$ respectively.\\
Now, let $(M,g)$ and $(N,h)$ be two closed Riemannian manifolds
and \bdism \varphi:(M,g)\longrightarrow(N,h) \edism an isometry,
we then have the following remarkable identity
\begin{align}\notag
\Delta_{g}\circ\varphi^{*}=\varphi^{*}\circ\Delta_{h}
\end{align}
see \cite{Gal}. Given a diffeomorphism $\psi:M^{n}\longrightarrow
M^{n}$ we have that \bdism \psi:(M^{n},\psi^{*}g)\longrightarrow
(M^{n},g) \edism is an isometry, we conclude that
$(M^{n},\psi^{*}g)$ and $(M^{n},g)$ have the same spectrum \bdism
Spec(g)=Spec(\psi^{*}g) \edism with eigenfunctions
$\left\{f_{k}\right\}$ and $\left\{\psi^{*}f_{k}\right\}$
respectively. We have thus that if $g(t)$ is a Ricci soliton on
$(M^{n},g_{0})$ then
\begin{align}\label{spec soliton}
Spec(g(t))=\frac{1}{\sigma(t)}Spec(g_{0}),
\end{align}
that is $Spec(g(t))$ is "proportional" to the initial spectrum
$Spec(g_{0})$ and shrinks, is stationary, or expands depending on
wheter $\epsilon$ is positive, zero, or negative. By \ref{spec
soliton}, we can also compute explicitely the change rate of each
eigenvalue of the Laplacian on a soliton \bdism
\frac{d\lambda}{dt}=-\frac{\sigma^{'}(t)}{\sigma(t)^{2}}=-\frac{2\epsilon}{(1+2\epsilon
t)^{2}}. \edism It could be now interesting to combine the above
formula with proposition \ref{evolution1} and try to study some
geometrical aspect of the soliton. This will be studied elsewhere.
\section{Remark}
The variation formula of proposition \ref{evolution1} also appears
in an equivalent form in \cite{Qian}, see section 7 theorem 7.3,
where it is derived as a consequence of a general formula for the
variation of the $L^{2}$-norm of a time dependent function on a
manifold evolving under the Ricci flow.

\section{Acknowledgements}
The author would like to thank Professor Claudio Procesi for
suggesting the study of the Ricci flow and Professor Stefano
Marchiafava for many valuable advices and suggestions. In
particular, Professor Stefano Marchiafava first pointed out to the
author the problem of the spectrum variation under a Ricci
deformation of the metric.

\end{document}